\documentclass[a4paper]{article}
\usepackage{ucs}
\usepackage[left=3cm,right=3cm,top=2cm,bottom=2cm]{geometry}
\usepackage{amsmath}
\usepackage{amsthm}
\usepackage{amssymb}
\usepackage{pgfplots}
\usepackage{pgfplotstable}
\usepackage{hyperref}
\usepackage{booktabs}
\usepackage{textcomp}
\usepackage{authblk}
\usepackage[sort]{cite}
\title{Adaptive time-stepping for aggregation-shattering kinetics}

\author[1,2]{Sergey A. Matveev\thanks{matseralex@cs.msu.ru}}
\author[1]{Viktor Zhilin}
\author[1]{Alexander P. Smirnov}
\affil[1]{Faculty of Computational Mathematics and Cybernetics, Lomonosov MSU, Moscow, Russia} 
\affil[2]{Marchuk Institute of Numerical Mathematics of the Russian Academy of Sciences, Moscow, Russia} 
\date{}
\begin{document}

\maketitle{}
\begin{abstract}
We propose an experimental study of adaptive time-stepping methods for efficient modeling of the aggregation-fragmentation kinetics. Precise modeling of this phenomena usually requires utilization of the large systems of nonlinear ordinary differential equations and intensive computations. We concentrate on performance of three explicit Runge-Kutta time-integration methods and provide simulations for two types of problems: finding of equilibrium solutions and simulations for kinetics with periodic solutions. The first class of problems may be analyzed through the relaxation of the solution to the stationary state after large time. In this case, the adaptive time-stepping may help to reach it using big steps reducing cost of the calculations without loss of accuracy. In the second case, the problem becomes numerically unstable at certain points of the phase space and may require tiny steps making the simulations very time-consuming. Adaptive criteria allows to increase the steps for most of points and speedup simulations significantly.
\end{abstract}


\maketitle

\section{Introduction}

Kinetic processes of aggregation and fragmentation often can be described with use of nonlinear ordinary differential equations. These processes are widespread in nature and important for many industries \cite{krapivsky2010kinetic, aloyan2014dynamics}. The smallest possible aggregates in such physical systems are usually called monomers. They may coalesce due to the collisions and generate larger polymer particles. These particles can also collide with each other further and grow up to larger particles consisting of thousands or even millions of monomers. If one knows or fixes the kinetic rates $K_{i,j}$ (kernel) for the reactions $[i]+[j] \rightarrow [i+j]$ with some formal expression (e.g. $K_{i,j}=(i/j)^{1/3}+(j/i)^{1/3}+2$) then the basic equations are well-known as Smoluchowski equations \cite{Smol17}:
\begin{equation}\label{eq:Smolbasic}
\frac{d n_s}{dt} = \underbrace{\frac{1}{2} \sum_{i+j=s} K_{i,j} n_i n_j}_{\text{"birth" of particle of size } s} - \underbrace{n_s \sum_{i =1}^{\infty} K_{s,i} n_i,}_{\text{"death" of particles of size } s}s = 1,2,\ldots,\infty.
\end{equation}
These equations describe evolution of the concentrations $n_s(t)$ of the particles of size $s$ due to their "birth" after the coalescence of smaller aggregates and due to their "death" after merging events. All interactions in such system are assumed to be pairwise, mass-conserving and spatially homogeneous.

Fixing some initial conditions as $n_s(t=0)$ one gets the Cauchy problem for this formally infinite system of nonlinear ODEs. This problem can be solved analytically only in the very rare special cases (e.g. constant kinetic rates with monodisperse or exponential initial conditions \cite{lushnikov1976evolution}). Hence, there is a certain need in efficient numerical algorithms for their investigation. Asymptotic and scaling features of the solutions $n_s(t)$ for $s\gg 1$, $t \gg 1$ are important for physicists \cite{leyvraz2003scaling, kang1986long} but their justification often may require utilization of enormous finite sub-systems
\begin{equation}\label{eq:SmolbasicTrunc}
\frac{d n_s}{dt} = \underbrace{\frac{1}{2} \sum_{i+j=s} K_{i,j} n_i n_j}_{\text{"birth" of particle of size } s} - \underbrace{n_s \sum_{i =1}^{M} K_{s,i} n_i,}_{\text{"death" of particles of size } s} s = 1,2,\ldots, M.
\end{equation}
These systems of ODEs with $M \gg 1$ allow to approximate the solution of initial equations \eqref{sec:exp} quite well for the large time segments \cite{kang1986long, mcleod1962infinite, kaushik2023steady}. Straight-forward computations with classical Runge-Kutta methods for these large systems take too much computing resources that should be likely reduced (each computation of the right-hand side takes $O(M^2)$ operations).  A family of efficient methods allowing to evaluate the right-hand side  $O(M R \log M)$ operations has been proposed for the case of the low-rank kernels 
$$K_{i,j} = \sum_{\alpha=1}^R U_{i, \alpha} V_{\alpha, j},R \ll M.$$
These methods allow to do computations with hundreds of thousands kinetic equations on basic laptops \cite{matveev2015fast}. Coupling of the low-rank decompositions with the time-integration methods has already been utilized for a broad class of problems. They include the irreversible aggregation with sources and sinks of particles \cite{matveev2020oscillating}
\begin{equation}\label{eq:SmolSources}
\frac{d n_s}{dt} = \frac{1}{2} \sum_{i+j=s} K_{i,j} n_i n_j - n_s \sum_{i =1}^{M} K_{s,i} n_i + \underbrace{P_k}_{\text{sources}}, \quad s = 1,2,\ldots, M,
\end{equation}
and aggregation-shattering kinetics \cite{matveev2017oscillations, brilliantov2015size}
\begin{eqnarray}
\notag
& \dfrac{d n_1}{dt} = \underbrace{- n_1 \sum\limits_{i=1}^{M} K_{1,i} n_i}_{\text{collusional aggregation}} +~ \underbrace{ \frac{\lambda}{2}\sum\limits_{i \geq 2}\sum\limits_{j\geq 2} (i+j) K_{i,j} n_i n_j + \lambda \sum\limits_{j\geq 2} j K_{1,j} n_j}_{\text{new monomers after collisional shattering events}}.\\
\label{eq:AggShatt} & ~\\ 
\notag
&\dfrac{d n_s}{dt} = \underbrace{\dfrac{1}{2} \sum\limits_{i+j=s} K_{i,j} n_i n_j- n_s \sum\limits_{i =1}^{M} K_{s,i} n_i}_{\text{collisional aggregation terms}} ~ - \underbrace{\lambda n_s \sum\limits_{i =1}^{M} K_{s,i}n_i, }_{\text{shattering into monomers}} \quad  s = 2,3,\ldots, M. 
\end{eqnarray}
Other applications also include temperature-dependent Smoluchoswki equations \cite{osinsky2020low, brilliantov2015size} and spatially inhomogeneous problems \cite{larchenko2023application}. 

Despite application of the low-rank decompositions leads to the drastic speedup of computations, most of interesting experiments have been done with the simplest constant time-steps that often might be a bottleneck for research. 

In case of the irreversible aggregation problem \eqref{eq:SmolSources} with multiple sources of particles \cite{matveev2020oscillating} there exist complex stationary particle size distributions that seem to be extremely hard for accurate construction with alternative approaches such as miscellaneous Monte Carlo methods \cite{smith1998constant, kruis2000direct, debry2003stochastic, sabelfeld1996stochastic} or coarse-graining approaches \cite{lee2000validity, stadnichuk2015smoluchowski}. Unfortunately, these stationary and quasi-stationary solutions may require to produce simulations for extremely large final times and number equations \cite{matveev2020oscillating, timokhin2019newton}.

At the same time, the dynamic oscillations are possible in aggregation models \cite{matveev2017oscillations, ball2012collective}. E. g. the simulations for the periodic solutions arising in the aggregation-shattering kinetics \eqref{eq:AggShatt} require very tiny time-steps \cite{matveev2017oscillations}. Otherwise, the simulations with constant time-steps become numerically unstable and crash. 

Adaptive schemes for the numerical integration of differential equations allow to automatically select the time-step depending on the characteristics of the system. These schemes are especially useful when the system has a variable sensitivity to changes over time or when an accurate solution is needed in certain time domains. One example of adaptive integration schemes is the Runge-Kutta method with automatic step selection (adaptive Runge-Kutta method see e.g. \cite{bakhvalov2000chislennye, driscoll2017fundamentals}). In this method, the integration step is varied in such a way as to control the approximation error. For small changes in the system, the integration step can be increased to save the computing resources, and for the fast dynamics the integration step can be reduced to ensure the accuracy of the solution.

Adaptive time integration schemes are widely used in various fields, such as modeling of dynamic systems, numerical solution of differential equations in physics, biology, economics and other sciences. They allow you to effectively and accurately simulate the behavior of systems with variable conditions.

There exists an old empirical criteria for stable calculations with Euler scheme for the choice of time-step
$$
\tau(t) \leq a \cdot \left[ \max  \sum_{j=1}^{M} K_{i,j} n_j(t) \right]^{-1}
$$
with $a \lesssim 1$. In practice, one may set $1/10<a <1/4$ but reaching the peak performance might require additional fine-tuning \cite{kang1986long, osinsky2020low}. Such a criteria is quite popular but it cannot be generalized for the higher order methods as well as has been elaborated only for the irreversible aggregation models.

We provide a numerical investigation showing that adaptive time-steps can be successfully applied to these problems. In this work  we
\begin{enumerate}
\item implement the adaptive time-steps for the Second, Fourth order Runge-Kutta and for the Runge-Kutta-Fehlberg methods in application to aggregation kinetics,
\item show speedup of calculations with adaptive-time steps for the problems with stationary and oscillating solutions by dozens of times,
\item present that the relaxation dynamics of the solutions to the stationary state for the large $T \gg 1$ can be tracked with adaptive-time steps.
\end{enumerate}

\section{Methods}\label{sec:methods}

In order to solve any of listed problems \eqref{eq:Smolbasic} -- \eqref{eq:AggShatt} with any time-integration method, we firstly re-write them in a compact operator form:
\begin{eqnarray}
\label{eq:vecform}
\left\lbrace \begin{matrix}
\dfrac{d\bf{n}}{dt} = \bf{S}(\bf{n})\\ 
~\\
\bf{n}(t=0) = \bf{n}_0
\end{matrix}\right. ,
\end{eqnarray}
where $\bf{n}(t)$ corresponds to the vector of concentrations $[n_1(t) , \ldots, n_M(t)]^{\top}$ and  $\bf{S}(\bf{n})$ denotes the right-hand side for the system of target ODEs. In this work, we do not discuss the complexity of evaluation of $\bf{S}(\bf{n})$ but only refer to the original work \cite{matveev2015fast}, where we have shown that it is $O(M R \log M)$ operations for the low-rank kernels with rank $R \ll M$. We utilize this approach and concentrate on investigation of performance of the time-integration methods. Namely, we utilize the following explicit schemes:
\begin{itemize}
\item second-order Runge-Kutta method (RK2)
\begin{align*}
& \bf{k}_1 = \tau \cdot \bf{S}(\bf{n}^k), \\
& \bf{k}_2 = \tau \cdot S(\bf{n}^{k}+\bf{k}_1), \\
& \bf{n}^{k+1} = \bf{n}^k + \dfrac{\bf{k}_1 + \bf{k}_2}{2}.
\end{align*}
\item fourth-order Runge-Kutta method (RK4)
\begin{align*}
& \bf{k}_1 = \tau \cdot \bf{S}(\bf{n}^{k})\\
& \bf{k}_2 = \tau \cdot \bf{S}(\bf{n}^{k} + \frac{1}{2} \bf{k}_1 ), \\
& \bf{k}_3 = \tau \cdot \bf{S}(\bf{n}^{k} + \frac{1}{2} \bf{k}_2 ), \\
& \bf{k}_4 = \tau \cdot \bf{S}(\bf{n}^{k} + k_3 ), \\
& \bf{n}^{k+1} = \bf{n}^{k} + \dfrac{\bf{k}_1 + 2 \bf{k}_2 + 2 \bf{k}_3 + \bf{k}_4}{6}.
\end{align*} 
\item Runge–Kutta–Fehlberg method (RKF45) \cite{fehlberg1970classical}
\begin{align*}
& \bf{k}_1 = \tau \cdot \bf{S}\left(n^{k}\right), \\
& \bf{k}_2 = \tau \cdot \bf{S}\left(n^{k} + \frac{1}{4}k_1\right), \\
& \bf{k}_3 = \tau \cdot \bf{S}\left(n^{k} + \frac{3 k_1 + 9 k_2}{32}\right), \\
& \bf{k}_4 = \tau \cdot \bf{S}\left(n^{k} + \frac{1932 k_1 - 7200 k_2 + 7296 k_3}{2197}\right), \\
& \bf{k}_5 = \tau \cdot \bf{S}\left(n^{k} + \frac{439}{216} k_1 - 8 k_2 + \frac{3680}{513} k_3 - \frac{845}{4104} k_4\right),\\
& \bf{k}_6 = \tau \cdot \bf{S}\left(n^{k} - \frac{8}{27} k_1 + 2 k_2 - \frac{3544}{2565} k_3 + \frac{1859}{4104} k_4 - \frac{11}{40} k_5\right),\\ 
& \bf{n}_{RK5} = \underbrace{\bf{n}^{k} + \frac{16}{135} \bf{k}_1 + \frac{6656}{12825} \bf{k}_3 + \frac{28561}{56430} \bf{k}_4 - \frac{9}{50} \bf{k}_5 + \frac{2}{55} \bf{k}_6}_{\text{step with the fifth order of accuracy for adaptive criteria. }} \\
& \bf{n}^{k+1} = \underbrace{\bf{n}^{k} + \frac{25}{216} \bf{k}_1 + \frac{1408}{2565} \bf{k}_3 + \frac{2197}{4104} \bf{k}_4 - \frac{1}{5} \bf{k}_5}_{\text{final step with the fouth order of accuracy, if criteria is fullfilled}}  \\
\end{align*}
\end{itemize} 
Thus, we seek to apply the adaptive criteria for selection of the time-integration step $\tau$ in order to reduce number of calls of $\bf{S}(\textbf{n})$ during the simulations. 

In case of the RK2 and the RK4 methods we utilize the general trick for automatic selection of the step-size (see e.g.  Ch.~8,~Par.~3 of \cite{bakhvalov2000chislennye} or \cite{driscoll2017fundamentals}) doing a pair of steps with $\tau/2$ obtaining a vector $\widehat{\bf{n}}^{k+1}$ and a single step $\tau$ getting $\bf{n}^{k+1}$. After it, the relative convergence error should be investigated according to basic theory and the time step is updated according to simple equation \cite{saleh2023adaptive}
\begin{equation}\label{eq:Apaptive_crit}
\tau_{new} := s \cdot \tau_{old} \cdot \left( \frac{tol}{||\bf{\widehat{n}^{k+1}- \bf{n}^{k+1}} ||_2}\right)^{1/p},
\end{equation}
where $tol$ is the tolerance specified by the user, $p$ is the order of the method ($p=2$ for the RK2 and $p=4$ for the RK4), $||\cdot||_2$ is Euclidian norm and  $0< s <1$ is a safety factor (we use $s = 1/4$).

Finally, if the error higher than pre-selected tolerance the time-step can become twice smaller and the refinement test has to be repeated. If the error is less than the tolerance criteria, then the next step can be done with twice larger $\tau$ and the refinement criteria should be checked again. 

For the RKF5 we follow the original criterion from paper by Felberg \cite{fehlberg1970classical}: the current iteration of the method has to be repeated until the estimate of the error becomes less than the maximum permissible as we set $tol$ parameter. In this case, on each attempt, the step decreases in proportion to the fifth root of the ratio of the maximal permissible error to the current one. For example, if the error exceeds the permissible by two times, then the step is reduced by approximately 20\%. If the error becomes less than the target level then the we increase the step size according to the rule. For example, if the error is two times less than the permissible, then the step increases by about 3\%. If the errors are equal, then the step is reduced by 10\%.

\section{Numerical experiments}\label{sec:exp}

\subsection{Complex steady-state for irreversible aggregation with multiple sources}

\begin{figure}[ht!]
    \centering
    \includegraphics[width=0.32\textwidth]{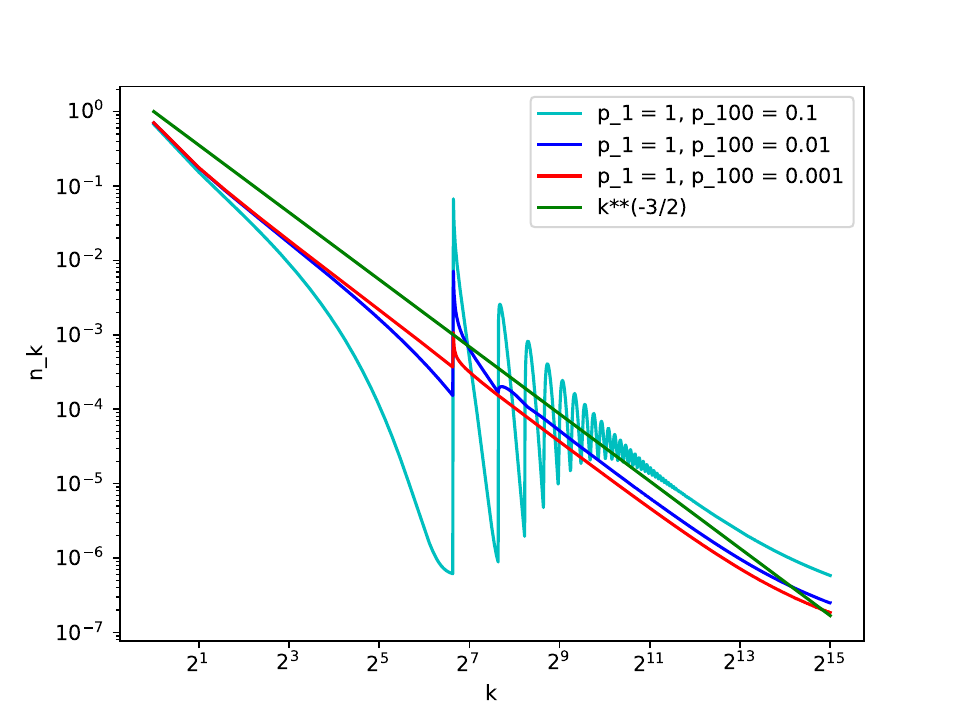}
    \includegraphics[width=0.32\textwidth]{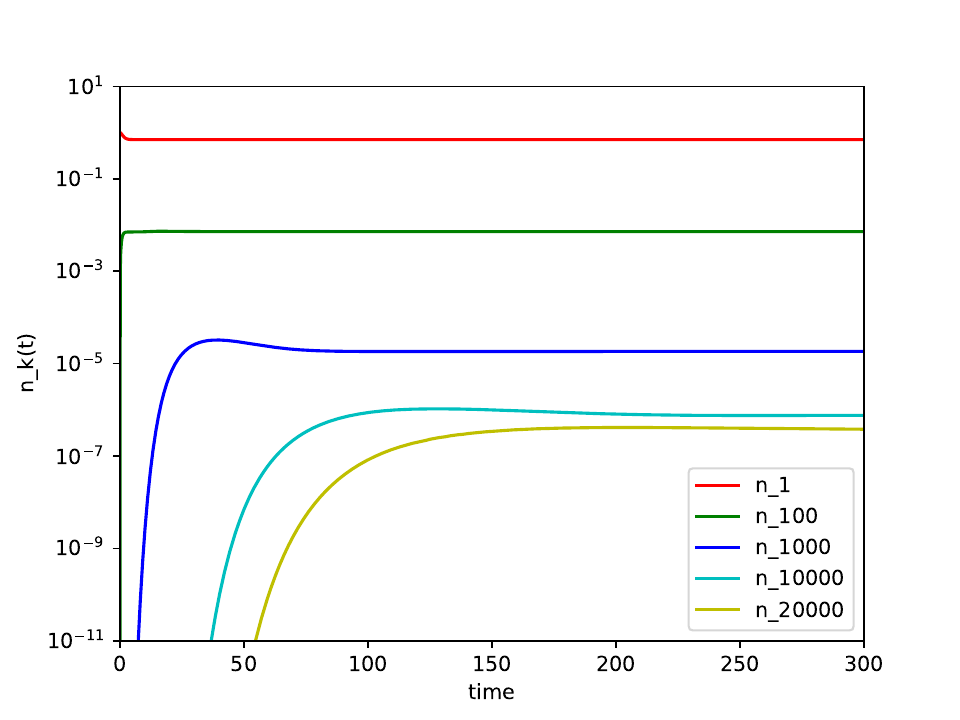}
    \includegraphics[width=0.32\textwidth]{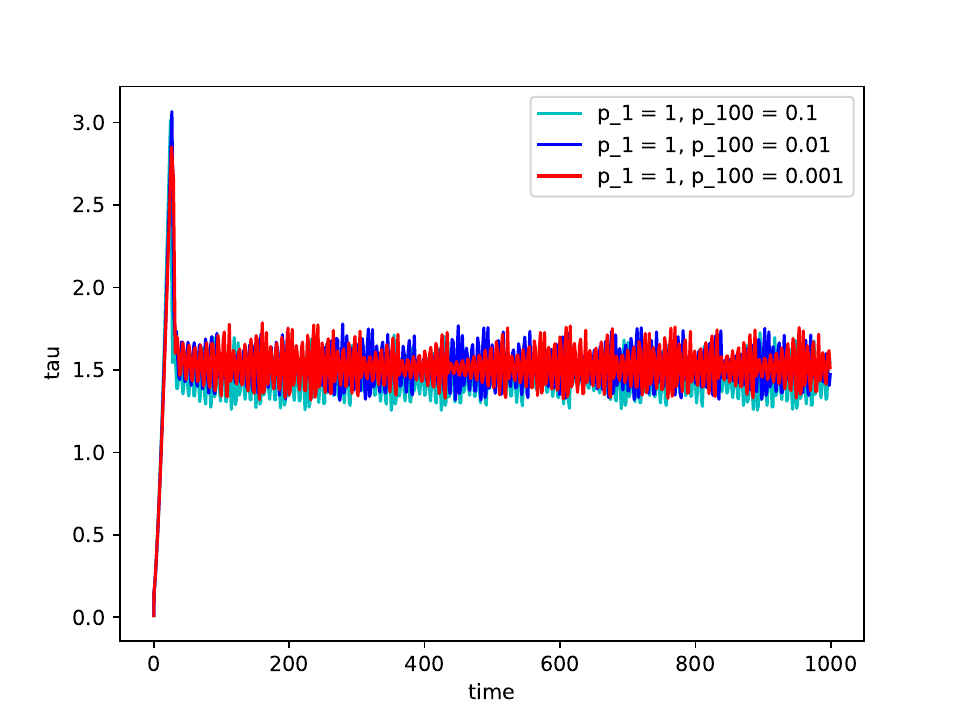}
    \caption{The problem with monodisperse conditions $M = 32768$ and a constant power of the source of monomers $P_1=1$ and several values of $P_{100}$ from $10^{-3}$ to $10^{-1}$. Left panel: the numerical stationary solutions oscillate for a wide range of the particle sizes in agreement with \cite{matveev2020oscillating} and relax to the scaling $n_k \simeq k^{-3/2}$, $k \gg 1$. Center panel: dynamics of the concentrations $n_k(t)$ for multiple masses. Right panel: the adaptive time-step increase quite rapidly and oscillate around unity for a long segment of model time.}
    \label{fig:two_sources_const_kernel}
\end{figure}

At first, we study the performance of the adaptive time-integration methods in application to irreversible aggregation equations with multiple sources \eqref{eq:SmolSources}. These equations can be studied with various kernels including the constant 
$$
K_{i,j} \equiv 1,
$$
free-molecular 
$$
K_{i,j} = (i^{1/3} + j^{1/3})^2 \sqrt{\frac{1}{i} + \frac{1}{j}}
$$
and family of generalized Brownian kernels
$$
K_{i,j} = \left(\frac{i}{j}\right)^{\alpha} + \left(\frac{j}{i}\right)^{\alpha}, \quad 0 \leq \alpha < 1/2, \quad c \geq 0.
$$
In recent work \cite{matveev2020oscillating}, we have found complex stationary solutions for these equations for all types of these kernels for a case of two constant sources of particles
$$
P_k = \left\lbrace \begin{matrix}
    1, & k = 1, \\
    p, & k = 100, \\
    0, & k \neq 1, 100.
\end{matrix} \right.
$$
Such complex stationary particle size distributions (see Fig. \ref{fig:two_sources_const_kernel}) correspond to the equilibrium state between the injection and sink processes for $t \rightarrow \infty$ leading to a system of nonlinear equations
\begin{equation}\label{eq:SmolSourcesSteady}
0 = \frac{1}{2} \sum_{i+j=s} K_{i,j} n_i n_j - n_s \sum_{i =1}^{M} K_{s,i} n_i + P_k, \quad s = 1,2,\ldots, M.
\end{equation}
The solutions of this system do not depend on initial conditions for ODEs (hence, we can them as monosdisperse $n_k(t=0)=\delta_{k,1}$, where $\delta_{i,j}$ is a Kronecker symbol). Finding such particle size distributions is an extremely tough problem via Monte Carlo simulations due to their highly oscillating form and long relaxation times. They have been found with use of deterministic iterative methods in application to \eqref{eq:SmolSourcesSteady}: one can apply Anderson acceleration \cite{matveev2020oscillating} or Newton method \cite{timokhin2019newton}. The iterative approach is very efficient and can be easily coupled with low-rank methods for evaluation of the right-hand size $\bf{S}(\bf(n))$. It allows to obtain the numerical solutions for the large systems even with million of equations.

\begin{table}[ht!]
\centering
\begin{tabular}{|c|c|c|c|c|}
\textbf{Scheme} & Const. $\tau=0.01$ & \textbf{$tol=10^{-4}$} & \textbf{$tol=10^{-6}$} & \textbf{$tol=10^{-8}$} \\
\hline
RK2  & 131.5 / 40000 & 2.737 / 784 & 17.66 / 4888 & 135.1 / 46882 \\ 
RK4  & 255.4 / 80000 & 3.633 / 1056 & 19.18 / 5780 & 177.5 / 54180 \\ 
RKF45 & 397.6 / 120000 & 3.223 / 876 & 3.557 / 972 & 4.579 / 1260 \\
\end{tabular}
\caption{Computational times in seconds and numbers of requests to the right-hand side for the constant kernel $K_{i,j}=1$,  $M = 32768$ for $T \in [0, 1000]$}
\label{tab:Sources_const}
\end{table}

\begin{table}[ht!]
\centering
\begin{tabular}{|c|c|c|c|c|}
\textbf{Scheme} & Const. $\tau=0.01$ & \textbf{$tol=10^{-4}$} & \textbf{$tol=10^{-6}$} & \textbf{$tol=10^{-8}$} \\
\hline
RK2  & 279.33 / 40000 & 35.52 / 5036 & 51.61 / 8784 & 305.94 / 49862 \\ 
RK4  & 482.12 / 80000 & 46.95 / 7508 & 70.60 / 11792 & 340.5 / 58788 \\ 
RKF45 & 690.5 / 120000 & 58.23 / 9966 & 57.25 / 9924 & 58.29 / 10056 \\
\end{tabular}
\caption{Computational times  in seconds and number of requests to the right-hand side for the Brownian kernel with $K_{i,j}=(i/j)^{1/3} + (j/i)^{1/3}$, $M = 32768$,  and $T \in [0, 1000]$}
\label{tab:Sources_Brownian}
\end{table}

However, the relaxation dynamics (see panel in the center of the Fig. \ref{fig:two_sources_const_kernel}) cannot be tracked with older methodology at all. At the same time, stationary solutions are reachable only for the large times and adaptive time-stepping methods become a necessary tool allowing to obtain the results in reasonable computing times. We sum up the results of our benchmarks in Table~\ref{tab:Sources_const} and Table~\ref{tab:Sources_Brownian}. Both computational times in seconds and number of evaluations of the right-hand side decrease significantly.

\subsection{Aggregation-shattering with dynamic oscillations}

The second example of a complex model for application of the adaptive time-stepping rules are the kinetic equations for aggregation-shattering processes \eqref{eq:AggShatt}.
Dynamic oscillations are possible for this model \cite{matveev2017oscillations} with non-local kernels
$$
K_{i,j} = \left(\frac{i}{j}\right)^{\alpha} + \left(\frac{j}{i}\right)^{\alpha}, \quad \alpha > 1/2,
$$
small shattering rates $0 < \lambda \leq \lambda^{*} \ll 1$ and various initial conditions (in this work we study the monodisperse case $n_k(t=0) = \delta_{k,1}$, where $\delta_{i,j}$ is a Kronecker symbol).

\begin{figure}[ht!]
    \centering
    \includegraphics[width=0.4\textwidth]{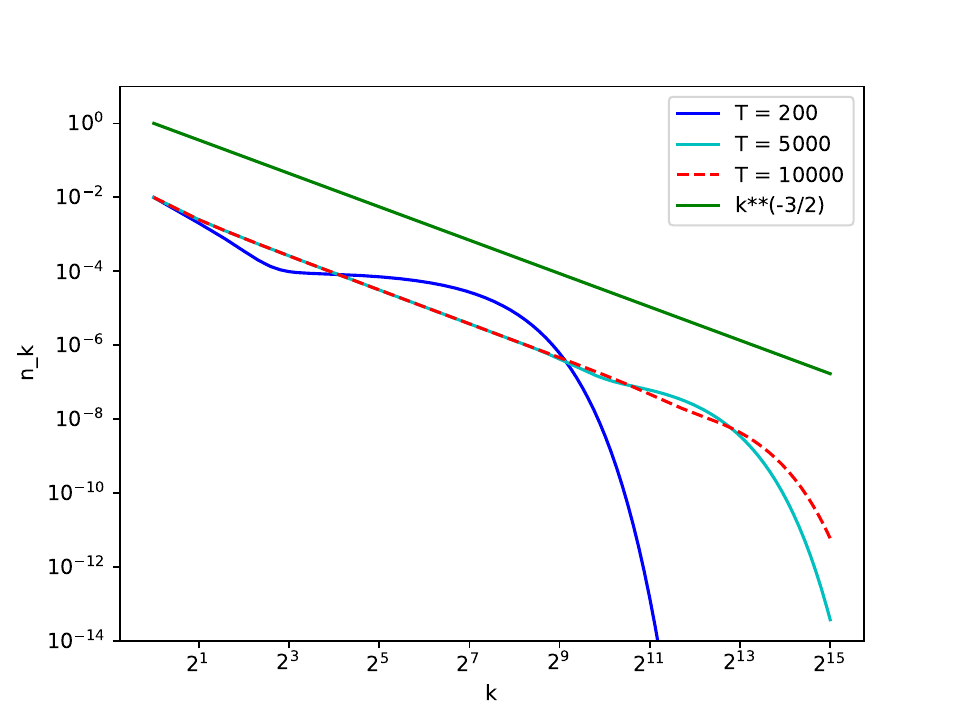}
    \includegraphics[width=0.4\textwidth]{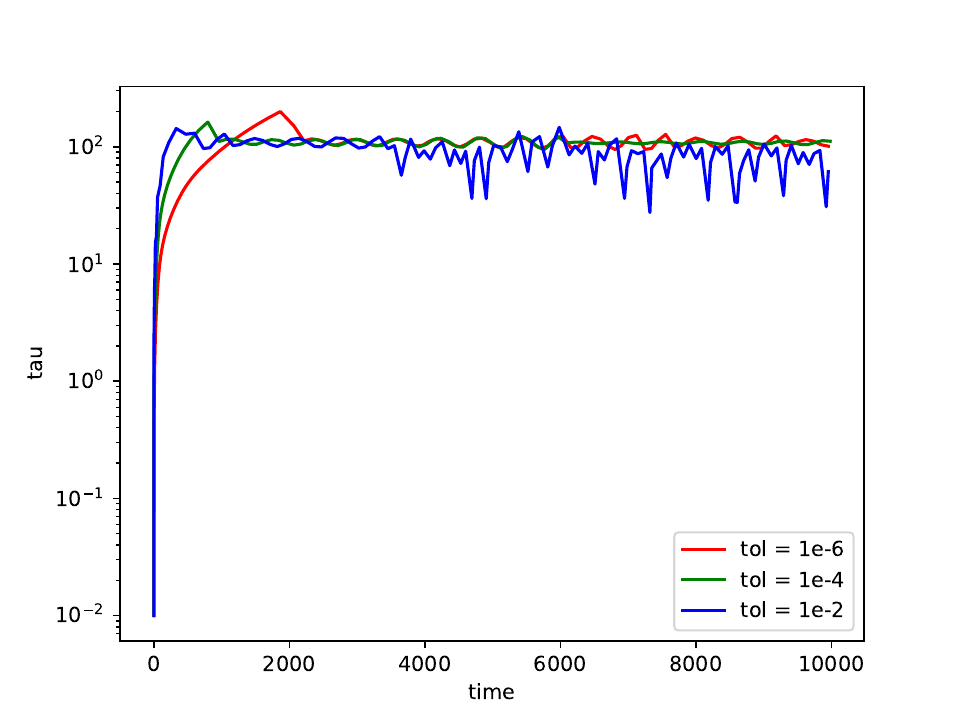}
    \caption{Numerical experiments for the constant kernel with $\alpha=0$, $M=32768$ and $\lambda= 0.01$ and monodisperse initial conditions. Left panel: very slow relaxation of the $n_k$ for $k\gg 1$ to the stationary distribution $n_k \propto k^{-3/2}$, the scaling becomes stable only for $T=10^4$. The solutions become close to the basic asymptotic form only for the very large $T$. Right panel: the time steps grow up to very large values if the adaptive time-stepping rule is applied.  }
    \label{fig:Aggr_shatt_steady}
\end{figure}

In Figure \ref{fig:Aggr_shatt_steady} we demonstrate the slow relaxation of the solution of the aggregation-shattering equations to the stationary state
$$
n_k \simeq k^{-\beta} \cdot e^{-\lambda^2 k}
$$
in agreement with previously elaborated theory \cite{brilliantov2015size} for the constant kernel ($\beta=-3/2$ for $\alpha=0$ and $\lambda=0.01$). We see that this relaxation is slow and takes a lot of time.  Utilization of the adaptive time-steps allows to track this transition dynamics with good accuracy within modest computing time. The previously elaborated iterative methods (see e.g. \cite{timokhin2019newton, stadnichuk2015smoluchowski}) allow to obtain only the numerical approximations of the final stationary particle size distributions in contrast with our approach. 

The small non-monotonicity in the numerical solution for $T=5000$ (see Fig. \ref{fig:Aggr_shatt_steady}) qualitatively agrees with the experimental data of the Voyager Radio Science Subsystem \cite{zebker1985saturn} for the Saturn's A ring, that can bee seen in the related plot for this model from \cite{brilliantov2015size}. Thus, we may infer that the real particle size distribution is close to the dynamic equilibrium predicted by theory but still needs much time in order to reach it.

\begin{figure}[ht!]
    \centering
    \includegraphics[width=0.4\textwidth]{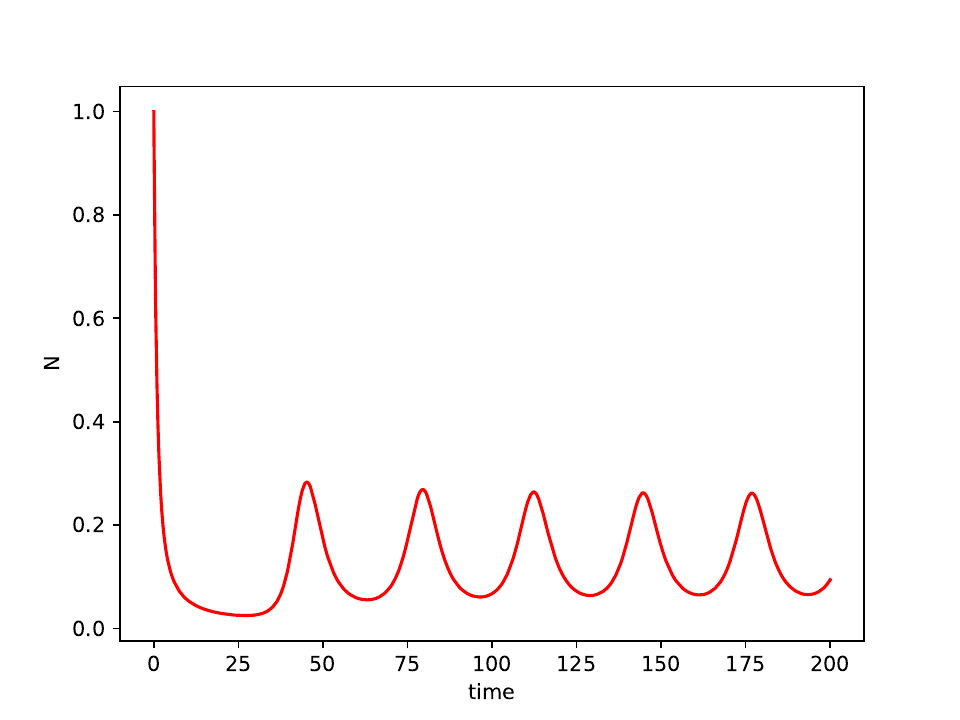}
    \includegraphics[width=0.4\textwidth]{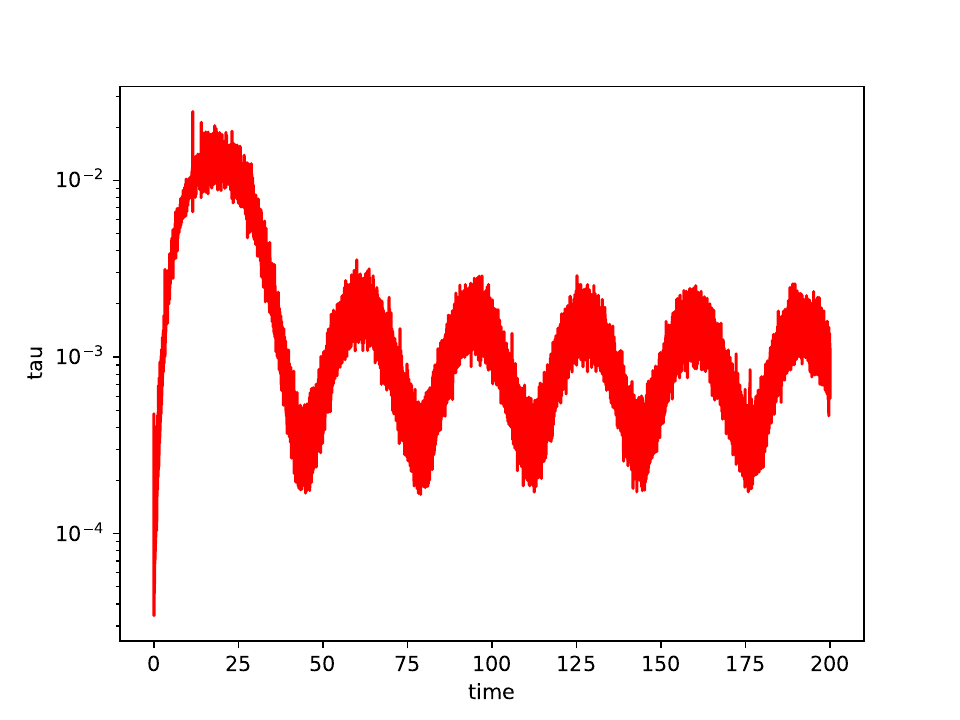}
    \caption{Numerical experiments for the generalized Brownian kernel with $\alpha=0.95$, $M=32768$ and $\lambda=0.01$ and monodisperse initial conditions. Left panel: oscillations of the total density $N(t)$ in agreement with baseline paper \cite{matveev2017oscillations}. Right panel: the adaptive time steps also oscillate with time and lead to speedup the computations.}
    \label{fig:Aggr_shatt_oscillations}
\end{figure}

In Figure \ref{fig:Aggr_shatt_oscillations} we demonstrate the results of our experiments for the kernel with $\alpha=0.98$ and $\lambda= 0.01$. In this case, the equations \eqref{eq:AggShatt} posses the dynamic oscillations. These oscillations have been found only numerically \cite{matveev2017oscillations} and seem to arise via Hopf bifurcation. The dynamic oscillations are also possible for the open irreversible aggregation model with the sources and sinks \cite{ball2012collective} but the aggregation-shattering kinetics is mass-conserving and such cyclic solutions even more interesting. 

The corresponding simulations basing on the second order Runge-Kutta method require tremendously small constant time-steps for stable calculations. Thus, even experiments for oscillations with twenty or forty thousand kinetic equations may take dozens of hours of computations.

On the right panel of Fig. \ref{fig:Aggr_shatt_oscillations} we see that the adaptive steps also oscillate with time and change almost by two orders of magnitude during each cycle. Application of the adaptive criteria allows to speedup these tough calculations by 10-15 times as we show in the Table \ref{tab:With_fragmentation}.  Basing on this observation we see that computations can be significantly accelerated without loss of their accuracy.

\begin{table}[ht!]
\centering
\begin{tabular}{|c|c|c|c|}
\textbf{Scheme} & \textbf{$tol=1e^{-4}$} & \textbf{$tol=1e^{-6}$} & \textbf{$tol=1e^{-8}$} \\
\hline
RK2 & 4160 & 4217 & 5250 \\
RK4 & 6399 & 6084 & 6520 \\
RKF45 & 9924 & 9700 & 9869 \\
\end{tabular}
\caption{Timings in seconds of numerical computations for the aggregation-shattering problem with dynamic oscillations ($\alpha=0.95$, $M=32768$, $\lambda=0.01$). The same calculations with constant time-steps $\tau = 5 \cdot 10^{-5}$ (the numerical integration becomes instable for the larger time-steps) with  requires approximately 23 hours with constant time-steps. }
\label{tab:With_fragmentation}
\end{table}

\section{Conclusion}

In this work we have studied the performance of several explicit time-integration methods with adaptive time-stepping criteria in application to problems of aggregation-fragmentation kinetics. Basing on our numerical experiments we obtain that complicated calculations for different models can be accelerated by dozens of times. 

We also demonstrate that application of the time-integration methods with adaptive time-steps allows to obtain the dynamics of relaxation process for the problems with stationary solutions using modest computing resources. Utilization of chosen dynamic time-stepping criteria allows us to exploit the higher-order Runge-Kutta methods instead of the well-known classical criteria for the simplest Euler scheme. 

Robust implementation of any implicit time-integration is an interesting direction for future research. It might be useful for the spatially inhomogeneous coagulation equations \cite{aloyan2014dynamics} playing important role in ecological modeling \cite{khodzher2021study}. This task seems to be challenging due to a certain need in efficient and highly accurate methods solving the systems of non-linear equations generated by Smoluchowski operator \cite{timokhin2019newton}. The cost of internal iterations within implicit time-steps must be smaller than straight-forward sequence of explicit time-steps. 

\section*{Funding and acknowledgements}
A. P. Smirnov was supported by Russian Science Foundation, grant project 24-11-00058 (\url{https://www.rscf.ru/project/24-11-00058/}). No funding was received to assist with the preparation of this manuscript by S. A. Matveev and Viktor Zhilin. S. A. Matveev acknowledges A. Osinsky for useful discussions.

\section*{Data availability and competing interests}
Data sharing not applicable to this article as no datasets were generated or analyzed during the current study. The source codes with implementation of the developed numerical methods are freely available upon request to the corresponding Author. The authors have no conflicts of interest to declare that are relevant to the content of this article.

\bibliography{bib}

\end{document}